\documentclass[12pt]{article}
\usepackage{amsfonts}
\usepackage{amsmath,amscd,amsthm,a4,amssymb}
\usepackage{amsmath,amscd,amsthm,a4,amssymb}

\newtheorem{theorem}{Theorem}
{}
\newtheorem{corollary}{Corollary}
{}

{}

{}

{}
\theoremstyle{plain}
{}

{}

{}

{}

{}

{}

{}

{}

{}

{}
\newtheorem{example}{Example}
{}

{}

{}

{}

{}

{}

\begin{document}
\begin{center}
{\Large \bf{General Rotational Surfaces with Pointwise 1- Type Gauss Map in
Pseudo- Euclidean Space E$_{2}^{4}$ \ }}
\end{center}
\centerline{\large Ferdag KAHRAMAN AKSOYAK  $^{1}$, Yusuf YAYLI $^{2}${\footnotetext{
{E-mail: $^{1}$ferda@erciyes.edu.tr(F. Kahraman Aksoyak ); $^{2}$yayli@science.ankara.edu.tr (Y.Yayli)}} }}

\

\centerline{\it $^{1}$Erciyes University, Department of Mathematics,
Kayseri, Turkey}
\centerline{\it $^{2}$Ankara University, Department of Mathematics,
Ankara, Turkey}

\begin{abstract}
In this paper, we study general rotational surfaces in the 4- dimensional
pseudo-Euclidean space $\mathbb{E}_{2}^{4}$ and obtain a characterization of
flat general rotation surfaces with pointwise 1-type Gauss map in $\mathbb{E}%
_{2}^{4}$ and give an example of such surfaces.
\end{abstract}

\begin{quote}\small
\end{quote}
\begin{quote}\small
\end{quote}

\section{Introduction}

A pseudo- Riemannian submanifold $M$ of the $m-$dimensional pseudo-Euclidean
space $\mathbb{E}_{s}^{m}$ is said to be of finite type if its position
vector $x$ can be expressed as a finite sum of eigenvectors of the Laplacian
$\Delta $ of $M$, that is, $x=x_{0}+x_{1}+...x_{k}$, where $x_{0}$ is a
constant map, $x_{1},...,x_{k}$ are non-constant maps such that $\Delta
x_{i}=\lambda _{i}x_{i},$ $\lambda _{i}\in $ $\mathbb{R}$, $i=1,2,...,k.$ If
$\lambda _{1},\lambda _{2},$...,$\lambda _{k}$ are all different, then $M$
is said to be of $k-$type. This definition was similarly extended to
differentiable maps in Euclidean and pseudo-Euclidean space, in particular,
to Gauss maps of submanifolds \cite{chen1}.

If a submanifold $M$ of a Euclidean space or pseudo-Euclidean space has
1-type Gauss map $G$, then $G$ satisfies $\Delta G=\lambda \left( G+C\right)
$ for some $\lambda \in \mathbb{R}$ and some constant vector $C.$ Chen and
Piccinni made a general study on compact submanifolds of Euclidean spaces
with finite type Gauss map and they proved that a compact hypersurface $M$
of $\mathbb{E}^{n+1}$ has 1-type Gauss map if and only if $M$ is a
hypersphere in $\mathbb{E}^{n+1}$ \cite{chen1}$.$

Hovewer the Laplacian of the Gauss map of several surfaces and hypersurfaces
such as a helicoids of the 1st,2nd and 3rdkind, conjugate Enneper's surface
of the second kind in 3- dimensional Minkowski space $E_{1}^{3}$,
generalized catenoids, spherical n-cones, hyperbolical n-cones and Enneper's
hypersurfaces in $E_{1}^{n}$ take the form namely,
\begin{equation}
\Delta G=f\left( G+C\right)
\end{equation}%
for some smooth function $f$ on $M$ and some constant vector $C.$ A
submanifold $M$ of a pseudo-Euclidean space $\mathbb{E}_{s}^{m}$ is said to
have pointwise 1-type Gauss map if its Gauss map satisfies $\left( 1\right) $
for some smooth function $f$ on $M$ and some constant vector $C.$ A
submanifold with pointwise 1-type Gauss map is said to be of the first kind
if the vector $C$ in $\left( 1\right) $ is zero\ vector. Otherwise, the
pointwise 1-type Gauss map is said to be of the second kind.

Surfaces in Euclidean space and in pseudo-Euclidean space with pointwise
1-type Gauss map were recently studied in \cite{chen}, \cite{choi1}, \cite%
{choi2}, \cite{choi3}, \cite{dursun2}, \cite{dursun3}, \cite{dursun4}, \cite%
{dursun5}, \cite{kim2}, \cite{niang1}, \cite{niang2}. Also Dursun and Turgay
in \cite{dursun1} gave all general rotational surfaces in $\mathbb{E}^{4}$
with proper pointwise 1-type Gauss map of the first kind and classified
minimal rotational surfaces with proper pointwise 1-type Gauss map of the
second kind. Arslan et al. in \cite{arslan1} investigated rotational
embedded surface with pointwise 1-type Gauss map. Arslan at el. in \cite%
{arslan2} gave necessary and sufficent conditions for Vranceanu rotation
surface to have pointwise 1-type Gauss map. Yoon in \cite{yoon2} showed that
flat Vranceanu rotation surface with pointwise 1-type Gauss map is a
Clifford torus and in \cite{yoon1} studied rotation surfaces in the
4-dimensional Euclidean space with finite type Gauss map. Kim and Yoon in
\cite{kim1} obtained the complete classification theorems for the flat
rotation surfaces with finite type Gauss map and pointwise 1-type Gauss map.
The authors in \cite{ak} studied flat general rotational surfaces in the 4-
dimensional Euclidean space $\mathbb{E}^{4}$ with pointwise 1-type Gauss map
and they showed that a non-planar flat general rotational surfaces with
pointwise 1-type Gauss map is a Lie group if and only if it is a Clifford
Torus.

In this paper, we study general rotational surfaces in the 4- dimensional
pseudo-Euclidean space $\mathbb{E}_{2}^{4}$ and obtain a characterization
for flat general rotation surfaces with pointwise 1-type Gauss map and give
an example of such surfaces.

\section{Preliminaries}

Let $E_{s}^{m}$ be the $m-$dimensional pseudo-Euclidean space with signature
$(s,m-s)$. Then the metric tensor $g$ in $E_{s}^{m}$ has the form
\begin{equation*}
g=\sum \limits_{i=1}^{m-s}\left( dx_{i}\right) ^{2}-\sum
\limits_{i=m-s+1}^{m}\left( dx_{i}\right) ^{2}
\end{equation*}%
where $(x_{1},...,x_{m})$ is a standard rectangular coordinate system in $%
E_{s}^{m}.$

Let $M$ be an $n-$dimensional pseudo-Riemannian submanifold of a $m-$%
dimensional pseudo-Euclidean space $\mathbb{E}_{s}^{m}.$ We denote
Levi-Civita connections of $\mathbb{E}_{s}^{m}$ and $M$ by $\tilde{\nabla}$
and $\nabla ,$ respectively. Let $e_{1},$...,$e_{n},e_{n+1},$...,$e_{m}$ be
an adapted local orthonormal frame in $\mathbb{E}_{s}^{m}$ such that $e_{1},$%
...,$e_{n}$ are tangent to $M$\ and $e_{n+1},$...,$e_{m}$ normal to $M.$ We
use the following convention on the ranges of indices: $1\leq i,j,k,$...$%
\leq n$, $n+1\leq r,s,t,$...$\leq m$, $1\leq A,B,C,$...$\leq m.$

Let $\omega _{A}$ be the dual-1 form of $e_{A}$ defined by $\omega
_{A}\left( X\right) =\left \langle e_{A},X\right \rangle $ and $\varepsilon
_{A}=\left \langle e_{A},e_{A}\right \rangle =\pm 1.$ Also, the connection
forms $\omega _{AB}$ are defined by%
\begin{equation*}
de_{A}=\sum \limits_{B}\varepsilon _{B}\omega _{AB}e_{B},\text{ \  \ }\omega
_{AB}+\omega _{BA}=0
\end{equation*}%
Then we have
\begin{equation}
\tilde{\nabla}_{e_{k}}^{e_{i}}=\sum \limits_{j=1}^{n}\varepsilon _{j}\omega
_{ij}\left( e_{k}\right) e_{j}+\sum \limits_{r=n+1}^{m}\varepsilon
_{r}h_{ik}^{r}e_{r}
\end{equation}%
and%
\begin{equation}
\tilde{\nabla}_{e_{k}}^{e_{s}}=-\sum \limits_{j=1}^{n}\varepsilon
_{j}h_{kj}^{s}e_{j}+\sum \limits_{r=n+1}^{m}\varepsilon _{r}\omega
_{sr}\left( e_{k}\right) e_{r},\text{....}D_{e_{k}}^{e_{s}}=\sum%
\limits_{r=n+1}^{m}\omega _{sr}\left( e_{k}\right) e_{r},
\end{equation}%
where $D$ is the normal connection, $h_{ik}^{r}$ the coefficients of the
second fundamental form $h.$

If we define a covariant differention $\tilde{\nabla}h$ of the second
fundamental form $h$ on the direct sum of the tangent bundle and the normal
bundle $TM\oplus T^{\perp }M$ of $M$ by
\begin{equation*}
\left( \tilde{\nabla}_{X}h\right) \left( Y,Z\right) =D_{X}h\left( Y,Z\right)
-h\left( \nabla _{X}Y,Z\right) -h\left( Y,\nabla _{X}Z\right)
\end{equation*}%
for any vector fields $X,$ $Y$ and $Z$ tangent to $M.$ Then we have the
Codazzi equation%
\begin{equation}
\left( \tilde{\nabla}_{X}h\right) \left( Y,Z\right) =\left( \tilde{\nabla}%
_{Y}h\right) \left( X,Z\right)
\end{equation}%
and the Gauss equation is given by
\begin{equation}
\left \langle R(X,Y)Z,W\right \rangle =\left \langle h\left( X,W\right)
,h\left( Y,Z\right) \right \rangle -\left \langle h\left( X,Z\right)
,h\left( Y,W\right) \right \rangle
\end{equation}%
where the vectors $X,$ $Y,$ $Z$ and $W$ are tangent to $M$ and $R$ is the
curvature tensor associated with $\nabla .$ The curvature tensor $R$
associated with $\nabla $ is defined by%
\begin{equation*}
R(X,Y)Z=\nabla _{X}\nabla _{Y}Z-\nabla _{Y}\nabla _{X}Z-\nabla _{\left[ X,Y%
\right] }Z.
\end{equation*}%
For any real function $f$ on $M$ the Laplacian $\Delta f$ of $f$ is given by
\begin{equation}
\Delta f=-\sum \limits_{i}\left( \tilde{\nabla}_{e_{i}}\tilde{\nabla}%
_{e_{i}}f-\tilde{\nabla}_{\nabla _{e_{i}}^{e_{i}}}f\right)
\end{equation}%
Let us now define the Gauss map $G$ of a submanifold $M$ into $G(n,m)$ in $%
\wedge ^{n}\mathbb{E}_{s}^{m},$ where $G(n,m)$ is the Grassmannian manifold
consisting of all oriented $n-$planes through the origin of $\mathbb{E}%
_{s}^{m}$ and $\wedge ^{n}\mathbb{E}_{s}^{m}$ is the vector space obtained
by the exterior product of $n$ vectors in $\mathbb{E}_{s}^{m}.$ Let $%
e_{i_{1}}\wedge ...\wedge e_{i_{n}}$ and $f_{j_{1}}\wedge ...\wedge
f_{j_{n}} $be two vectors of $\wedge ^{n}\mathbb{E}_{s}^{m},$ where $%
\left
\{ e_{1},\text{...,}e_{m}\right \} $ and $\left \{ f_{1},\text{...,}%
f_{m}\right \} $ are orthonormal bases of $\mathbb{E}_{s}^{m}$. Define an
indefinite inner product $\left \langle ,\right \rangle $ on $\wedge ^{n}%
\mathbb{E}_{s}^{m}$ by%
\begin{equation*}
\left \langle e_{i_{1}}\wedge ...\wedge e_{i_{n}},f_{j_{1}}\wedge ...\wedge
f_{j_{n}}\right \rangle =\det \left( \left \langle e_{i_{l}},f_{j_{k}}\right
\rangle \right) .
\end{equation*}%
Therefore, for some positive integer $t,$ we may identify $\wedge ^{n}%
\mathbb{E}_{s}^{m}$ with some Euclidean space $\mathbb{E}_{t}^{N}$ where $%
N=\left(
\begin{array}{c}
m \\
n%
\end{array}%
\right) .$ The map $G:M\rightarrow G(n,m)\subset E_{k}^{N}$ defined by $%
G(p)=\left( e_{1}\wedge ...\wedge e_{n}\right) \left( p\right) $ is called
the Gauss map of $M,$ that is, a smooth map which carries a point $p$ in $M$
into the oriented $n-$plane in $\mathbb{E}_{s}^{m}$ obtained from parallel
translation of the tangent space of $M$ at $p$ in $\mathbb{E}_{s}^{m}.$

\section{Flat Rotation Surfaces with Pointwise 1-Type Gauss Map in $%
E_{2}^{4} $}

In this section, we study the flat rotation surfaces with pointwise 1-type
Gauss map in the 4-dimensional pseudo-Euclidean space $E_{2}^{4}$. Let $%
M_{1} $ and $M_{2}$ be the rotation surfaces in $E_{2}^{4}$ defined by%
\begin{equation*}
\varphi \left( t,s\right) =%
\begin{pmatrix}
\cosh t & 0 & 0 & \sinh t \\
0 & \cosh t & \sinh t & 0 \\
0 & \sinh t & \cosh t & 0 \\
\sinh t & 0 & 0 & \cosh t%
\end{pmatrix}%
\left(
\begin{array}{c}
0 \\
x(s) \\
0 \\
y(s)%
\end{array}%
\right) ,
\end{equation*}%
\begin{equation}
M_{1}:\text{ }\varphi \left( t,s\right) =\left( y(s)\sinh t,x(s)\cosh
t,x(s)\sinh t,y(s)\cosh t\right) \text{\ }
\end{equation}%
and
\begin{equation*}
\varphi \left( t,s\right) =%
\begin{pmatrix}
\cos t & -\sin t & 0 & 0 \\
\sin t & \cos t & 0 & 0 \\
0 & 0 & \cos t & -\sin t \\
0 & 0 & \sin t & \cos t%
\end{pmatrix}%
\left(
\begin{array}{c}
x(s) \\
0 \\
y(s) \\
0%
\end{array}%
\right)
\end{equation*}%
\begin{equation}
M_{2}:\text{ }\varphi \left( t,s\right) =\left( x(s)\cos t,x(s)\sin
t,y(s)\cos t,y(s)\sin t\right)
\end{equation}%
where the profile curve of $M_{1}$ (resp. the profile curve of $M_{2})$ is
unit speed curve, that is, $\left( x^{\prime }(s)\right) ^{2}-\left(
y^{\prime }(s)\right) ^{2}=1.$ We choose a moving frame $%
e_{1},e_{2},e_{3},e_{4}$ such that $e_{1},e_{2}$ are tangent to $M_{1}$ and $%
e_{3},e_{4}$ are normal to $M_{1}$ and choose a moving frame $\bar{e}_{1},%
\bar{e}_{2},\bar{e}_{3},\bar{e}_{4}$ such that $\bar{e}_{1},\bar{e}_{2}$ are
tangent to $M_{2}$ and $\bar{e}_{3},\bar{e}_{4}$ are normal to $M_{2}$ which
are given by the following:
\begin{eqnarray*}
e_{1} &=&\frac{1}{\sqrt{\varepsilon _{1}\left( y^{2}(s)-x^{2}\left( s\right)
\right) }}\left( y\left( s\right) \cosh t,x\left( s\right) \sinh t,x(s)\cosh
t,y(s)\sinh t\right) \\
e_{2} &=&\left( y^{\prime }\left( s\right) \sinh t,,x^{\prime }\left(
s\right) \cosh t,x^{\prime }(s)\sinh t,,y^{\prime }(s)\cosh t\right) \\
e_{3} &=&\left( x^{\prime }(s)\sinh t,y^{\prime }(s)\cosh t,y^{\prime
}\left( s\right) \sinh t,x^{\prime }\left( s\right) \cosh t\right) \\
e_{4} &=&\frac{1}{\sqrt{\varepsilon _{1}\left( y^{2}(s)-x^{2}\left( s\right)
\right) }}\left( x(s)\cosh t,y(s)\sinh t,y\left( s\right) \cosh t,x\left(
s\right) \sinh t\right)
\end{eqnarray*}%
and
\begin{eqnarray*}
\bar{e}_{1} &=&\frac{1}{\sqrt{\varepsilon _{1}\left( y^{2}(s)-x^{2}\left(
s\right) \right) }}\left( -x\left( s\right) \sin t,x\left( s\right) \cos
t,-y(s)\sin t,y(s)\cos t\right) \\
\bar{e}_{2} &=&\left( x^{\prime }\left( s\right) \cos t,,x^{\prime }\left(
s\right) \sin t,y^{\prime }(s)\cos t,y^{\prime }(s)\sin t\right) \\
\bar{e}_{3} &=&\left( y^{\prime }\left( s\right) \cos t,y^{\prime }\left(
s\right) \sin t,x^{\prime }(s)\cos t,x^{\prime }(s)\sin t\right) \\
\bar{e}_{4} &=&\frac{1}{\sqrt{\varepsilon _{1}\left( y^{2}(s)-x^{2}\left(
s\right) \right) }}\left( y\left( s\right) \sin t,-y\left( s\right) \cos
t,x(s)\sin t,-x(s)\cos t\right)
\end{eqnarray*}%
where $\varepsilon _{1}\left( y^{2}(s)-x^{2}\left( s\right) \right) >0,$ $%
\varepsilon _{1}=\pm 1.$ Then it is easily seen that
\begin{eqnarray*}
\left \langle e_{1},e_{1}\right \rangle &=&-\left \langle e_{4},e_{4}\right
\rangle =\varepsilon _{1},\text{ \  \  \ }\left \langle e_{2},e_{2}\right
\rangle =-\left \langle e_{3},e_{3}\right \rangle =1 \\
-\left \langle \bar{e}_{1},\bar{e}_{1}\right \rangle &=&\left \langle \bar{e}%
_{4},\bar{e}_{4}\right \rangle =\varepsilon _{1},\text{ \  \  \  \  \ }\left
\langle \bar{e}_{2},\bar{e}_{2}\right \rangle =-\left \langle \bar{e}_{3},%
\bar{e}_{3}\right \rangle =1
\end{eqnarray*}%
we have the dual 1-forms as:
\begin{equation}
\omega _{1}=\varepsilon _{1}\sqrt{\varepsilon _{1}\left(
y^{2}(s)-x^{2}\left( s\right) \right) }dt\text{ \  \  \  \ and \  \  \  \ }\omega
_{2}=ds
\end{equation}%
and%
\begin{equation}
\bar{\omega}_{1}=-\varepsilon _{1}\sqrt{\varepsilon _{1}\left(
y^{2}(s)-x^{2}\left( s\right) \right) }dt\text{ \  \  \  \ and \  \  \  \ }\bar{%
\omega}_{2}=ds
\end{equation}%
By a direct computation we have components of the second fundamental form
and the connection forms as:%
\begin{eqnarray}
h_{11}^{3} &=&b(s),\ h_{12}^{3}=0,\ h_{22}^{3}=c(s) \\
h_{11}^{4} &=&0,\text{ \ }h_{12}^{4}=b(s),\text{ \ }h_{22}^{4}=0  \notag
\end{eqnarray}%
\begin{eqnarray}
\bar{h}_{11}^{3} &=&-b(s),\  \bar{h}_{12}^{3}=0,\  \bar{h}_{22}^{3}=c(s) \\
\bar{h}_{11}^{4} &=&0,\text{ \ }\bar{h}_{12}^{4}=b(s),\text{ \ }\bar{h}%
_{22}^{4}=0  \notag
\end{eqnarray}%
\begin{eqnarray}
\omega _{12} &=&\varepsilon _{1}a(s)\omega _{1},\text{ \  \ }\omega
_{13}=\varepsilon _{1}b(s)\omega _{1},\text{ \  \ }\omega _{14}=b(s)\omega
_{2} \\
\omega _{23} &=&c(s)\omega _{2},\text{ \  \ }\omega _{24}=\varepsilon
_{1}b(s)\omega _{1},\text{ \  \ }\omega _{34}=\varepsilon _{1}a(s)\omega _{1}
\notag
\end{eqnarray}%
\begin{eqnarray}
\bar{\omega}_{12} &=&\varepsilon _{1}a(s)\bar{\omega}_{1},\text{ \  \ }\bar{%
\omega}_{13}=\varepsilon _{1}b(s)\bar{\omega}_{1},\text{ \  \ }\bar{\omega}%
_{14}=b(s)\bar{\omega}_{2} \\
\bar{\omega}_{23} &=&c(s)\bar{\omega}_{2},\text{ \  \ }\bar{\omega}%
_{24}=-\varepsilon _{1}b(s)\bar{\omega}_{1},\text{ \  \ }\bar{\omega}%
_{34}=-\varepsilon _{1}a(s)\bar{\omega}_{1}  \notag
\end{eqnarray}%
By covariant differentiation with respect to $e_{1}$ and $e_{2}$ (resp. $%
\bar{e}_{1}$ and $\bar{e}_{2})$ a straightforward calculation gives:
\begin{eqnarray}
\tilde{\nabla}_{e_{1}}e_{1} &=&a(s)e_{2}-b(s)e_{3} \\
\tilde{\nabla}_{e_{2}}e_{1} &=&-\varepsilon _{1}b(s)e_{4}  \notag \\
\tilde{\nabla}_{e_{1}}e_{2} &=&-\varepsilon _{1}a(s)e_{1}-\varepsilon
_{1}b(s)e_{4}  \notag \\
\tilde{\nabla}_{e_{2}}e_{2} &=&-c(s)e_{3}  \notag \\
\tilde{\nabla}_{e_{1}}e_{3} &=&-\varepsilon _{1}b(s)e_{1}-\varepsilon
_{1}a(s)e_{4}  \notag \\
\tilde{\nabla}_{e_{2}}e_{3} &=&-c(s)e_{2}  \notag \\
\tilde{\nabla}_{e_{1}}e_{4} &=&-b(s)e_{2}+a(s)e_{3}  \notag \\
\tilde{\nabla}_{e_{2}}e_{4} &=&-\varepsilon _{1}b(s)e_{1}  \notag
\end{eqnarray}%
and
\begin{eqnarray}
\tilde{\nabla}_{\bar{e}_{1}}\bar{e}_{1} &=&-a(s)\bar{e}_{2}+b(s)\bar{e}_{3}
\\
\tilde{\nabla}_{\bar{e}_{2}}\bar{e}_{1} &=&\varepsilon _{1}b(s)\bar{e}_{4}
\notag \\
\tilde{\nabla}_{\bar{e}_{1}}\bar{e}_{2} &=&-\varepsilon _{1}a(s)\bar{e}%
_{1}+\varepsilon _{1}b(s)\bar{e}_{4}  \notag \\
\tilde{\nabla}_{\bar{e}_{2}}\bar{e}_{2} &=&-c(s)\bar{e}_{3}  \notag \\
\tilde{\nabla}_{\bar{e}_{1}}e_{3} &=&-\varepsilon _{1}b(s)\bar{e}%
_{1}+\varepsilon _{1}a(s)\bar{e}_{4}  \notag \\
\tilde{\nabla}_{\bar{e}_{2}}e_{3} &=&-c(s)\bar{e}_{2}  \notag \\
\tilde{\nabla}_{\bar{e}_{1}}e_{4} &=&-b(s)\bar{e}_{2}+a(s)\bar{e}_{3}  \notag
\\
\tilde{\nabla}_{\bar{e}_{2}}e_{4} &=&\varepsilon _{1}b(s)\bar{e}_{1}  \notag
\end{eqnarray}%
where
\begin{equation}
a(s)=\frac{x(s)x^{\prime }(s)-y(s)y^{\prime }(s)}{\varepsilon _{1}\left(
y^{2}(s)-x^{2}\left( s\right) \right) }
\end{equation}%
\begin{equation}
\text{\ }b(s)=\frac{x(s)y^{\prime }(s)-x^{\prime }(s)y(s)}{\varepsilon
_{1}\left( y^{2}(s)-x^{2}\left( s\right) \right) }
\end{equation}%
\begin{equation}
c(s)=x^{\prime \prime }(s)y^{\prime }(s)-x^{\prime }(s)y^{\prime \prime }(s)
\end{equation}%
The Gaussian curvature $K$ of $M_{1}$ and $\bar{K}$ that of $M_{2}$ are
respectively given by
\begin{equation}
K=\varepsilon _{1}b^{2}(s)-b(s)c(s)
\end{equation}%
and
\begin{equation}
\bar{K}=b(s)c(s)-\varepsilon _{1}b^{2}(s)
\end{equation}%
If the surfaces $M_{1}$ or $M_{2}$ is flat, then $(20)$ and $\left(
21\right) $ imply
\begin{equation}
b(s)c(s)-b^{2}(s)=0.
\end{equation}%
Furthermore, after some computations we obtain Gauss and Codazzi equations
for both surfaces $M_{1}$ and $M_{2}$
\begin{equation}
\varepsilon _{1}a^{2}\left( s\right) -a^{\prime }\left( s\right)
=b(s)c(s)-\varepsilon _{1}b^{2}(s)
\end{equation}%
and
\begin{equation}
b^{\prime }\left( s\right) =2\varepsilon _{1}a(s)b(s)-a(s)c(s)
\end{equation}%
respectively.

By using $\left( 6\right) ,$ $\left( 15\right) ,$ $\left( 16\right) $ and
straight-forward computation, the Laplacians $\Delta G$ and $\Delta \bar{G}$
of the Gauss map $G$ and $\bar{G}$ can be expressed as%
\begin{eqnarray}
\Delta G &=&-\left( 3b^{2}\left( s\right) +c^{2}\left( s\right) \right)
\left( e_{1}\wedge e_{2}\right) +\left( 2a(s)b(s)-\varepsilon
_{1}a(s)c(s)+c^{\prime }\left( s\right) \right) \left( e_{1}\wedge
e_{3}\right)  \notag \\
&&+\left( 3a(s)b(s)-\varepsilon _{1}b^{\prime }(s)\right) \left( e_{2}\wedge
e_{4}\right) +2\left( \varepsilon _{1}b(s)c(s)-b^{2}(s)\right) \left(
e_{3}\wedge e_{4}\right)
\end{eqnarray}%
\begin{eqnarray}
\Delta \bar{G} &=&-\left( 3b^{2}\left( s\right) +c^{2}\left( s\right)
\right) \left( e_{1}\wedge e_{2}\right) +\left( 2a(s)b(s)-\varepsilon
_{1}a(s)c(s)+c^{\prime }\left( s\right) \right) \left( e_{1}\wedge
e_{3}\right)  \notag \\
&&+\left( -3a(s)b(s)+\varepsilon _{1}b^{\prime }(s)\right) \left(
e_{2}\wedge e_{4}\right) +2\left( b^{2}(s)-\varepsilon _{1}b(s)c(s)\right)
\left( e_{3}\wedge e_{4}\right)
\end{eqnarray}

Now we investigate the flat rotation surfaces in $E_{2}^{4}$ with the
pointwise 1-type Gauss map satisfying $(1)$.

Suppose that the rotation surface $M_{1}$ given by the parametrization $(7)$
is a flat rotation surface. From $(20)$, we obtain that $b(s)=0$ or $%
\varepsilon _{1}b(s)-c(s)=0.$ We assume that $\varepsilon _{1}b(s)-c(s)\neq
0.$ Then $b(s)$ is equal to zero and $(24)$ implies that $a(s)c(s)=0.$ Since
$\varepsilon _{1}b(s)-c(s)\neq 0,$ it implies that $c(s)$ is not equal to
zero. Then we obtain as $a(s)=0.$ In that case, by using $(17)$ and $(18)$
we obtain that $\alpha \left( s\right) =\left( 0,x\left( s\right)
,0,y(s)\right) $ is a constant vector. This is a contradiction. Therefore $%
\varepsilon _{1}b(s)=c(s)$ for all $s.$ From $(14)$, we get
\begin{equation}
\varepsilon _{1}a^{2}\left( s\right) -a^{\prime }\left( s\right) =0
\end{equation}%
whose the trivial solution and non-trivial solution%
\begin{equation*}
a(s)=0
\end{equation*}%
and%
\begin{equation*}
a(s)=\frac{1}{-\varepsilon _{1}s+c},
\end{equation*}%
respectively. We assume that $a(s)=0.$ By $(24)$ $b=b_{0}$ is a constant and
$c=\varepsilon _{1}b_{0}$. In that case by using $(17)$, $(18)$ and $(19)$, $%
x$ and $y$ satisfy the following differential equations
\begin{equation}
x^{2}\left( s\right) -y^{2}(s)=\mu \text{ \  \ }\mu \text{ is a constant,}
\end{equation}%
\begin{equation}
x(s)y^{\prime }(s)-x^{\prime }(s)y(s)=-\varepsilon _{1}b_{0}\mu ,
\end{equation}%
\begin{equation}
x^{\prime \prime }y^{\prime }(s)-x^{\prime }(s)y^{\prime \prime
}=\varepsilon _{1}b_{0}.
\end{equation}%
From $(28)$ we may put
\begin{equation}
x\left( s\right) =\frac{1}{2}\varepsilon \left( \mu _{2}e^{\theta \left(
s\right) }+\mu _{1}e^{-\theta \left( s\right) }\right) ,\text{ \  \ }y\left(
s\right) =\frac{1}{2}\varepsilon \left( \mu _{2}e^{\theta \left( s\right)
}-\mu _{1}e^{-\theta \left( s\right) }\right) ,
\end{equation}%
where $\theta \left( s\right) $ is some smooth function, $\varepsilon =\pm 1$
and $\mu =\mu _{1}\mu _{2}$. Differentiating $(31)$ with respect to $s,$ we
have%
\begin{equation}
x^{\prime }(s)=\theta ^{\prime }(s)y(s),\text{ \ }y^{\prime }(s)=\theta
^{\prime }(s)x(s)
\end{equation}%
By substituting $(31)$ and $(32)$ into $(19)$, we get
\begin{equation*}
\theta \left( s\right) =-\varepsilon _{1}b_{0}s+d\text{, \  \ }d=const.
\end{equation*}%
And since the curve $\alpha $ is a unit speed curve, we have
\begin{equation*}
b_{0}^{2}\mu =-1.
\end{equation*}%
Since $\mu =-\frac{1}{b_{0}^{2}},$ $y^{2}(s)-x^{2}\left( s\right) >0.$ In
that case we obtain that $\varepsilon _{1}=1.$ Then we can write components
of the curve $\alpha $ as:
\begin{eqnarray}
x\left( s\right)  &=&\frac{1}{2}\varepsilon \left( \mu _{2}e^{\left(
-b_{0}s+d\right) }+\mu _{1}e^{-(-b_{0}s+d)}\right) ,\text{ } \\
y\left( s\right)  &=&\frac{1}{2}\varepsilon \left( \mu _{2}e^{\left(
-b_{0}s+d\right) }-\mu _{1}e^{-\left( -b_{0}s+d\right) }\right) ,\text{ }\mu
_{1}\mu _{2}=-\frac{1}{b_{0}^{2}}  \notag
\end{eqnarray}%
On the other hand, by using $(25)$ we can rewrite the Laplacian of the Gauss
map $G$ with $a(s)=0$ and $b=c=b_{0}$ as follows:%
\begin{equation*}
\Delta G=-4b_{0}^{2}\left( e_{1}\wedge e_{2}\right)
\end{equation*}%
that is, the flat surface $M$ is pointwise 1-type Gauss map with the
function $f=4b_{0}^{2}$ and $C=0.$ Even if it is a pointwise 1-type Gauss
map of the first kind.

Now we assume that $a(s)=\frac{1}{-\varepsilon _{1}s+c}.$ By using $%
c(s)=\varepsilon _{1}b(s)$ and $(24)$ we get
\begin{equation}
b^{\prime }\left( s\right) =\varepsilon _{1}a(s)b(s)
\end{equation}%
or we can write%
\begin{equation*}
\frac{b^{\prime }\left( s\right) }{b(s)}=\frac{\varepsilon _{1}}{%
-\varepsilon _{1}s+c},
\end{equation*}%
whose the solution
\begin{equation}
b(s)=\frac{\lambda }{\left \vert -\varepsilon _{1}s+c\right \vert },\text{ \  \
\ }\lambda \text{ is a constant.}
\end{equation}%
By using $(25)$ we can rewrite the Laplacian of the Gauss map $G$ with the
equations $c(s)=\varepsilon _{1}b(s),$ $b^{\prime }\left( s\right)
=\varepsilon _{1}a(s)b(s)$ and $a^{\prime }\left( s\right) =\varepsilon
_{1}a^{2}\left( s\right) $%
\begin{equation}
\Delta G=-4b^{2}\left( s\right) \left( e_{1}\wedge e_{2}\right)
+2a(s)b(s)\left( e_{1}\wedge e_{3}\right) +2a(s)b(s)\left( e_{2}\wedge
e_{4}\right) .
\end{equation}%
We suppose that the flat rotational surface $M_{1}$ has pointwise 1-type
Gauss map. From $(1)$ and $(36)$, we get%
\begin{equation}
-4\varepsilon _{1}b^{2}\left( s\right) =f\varepsilon _{1}+f\left \langle
C,e_{1}\wedge e_{2}\right \rangle
\end{equation}%
\begin{equation}
-2\varepsilon _{1}a(s)b(s)=f\left \langle C,e_{1}\wedge e_{3}\right \rangle
\end{equation}%
\begin{equation}
-2\varepsilon _{1}a(s)b(s)=f\left \langle C,e_{2}\wedge e_{4}\right \rangle
\end{equation}%
Then, we have%
\begin{equation}
\left \langle C,e_{1}\wedge e_{4}\right \rangle =0,\text{ }\left \langle
C,e_{2}\wedge e_{3}\right \rangle =0,\text{ }\left \langle C,e_{3}\wedge
e_{4}\right \rangle =0
\end{equation}%
By using $(38)$ and $(39)$ we obtain
\begin{equation}
\left \langle C,e_{1}\wedge e_{3}\right \rangle =\left \langle C,e_{2}\wedge
e_{4}\right \rangle
\end{equation}%
By differentiating the first equation in $(41)$ with respect to $e_{1}$ and
by using the third equation in $(41)$ and $(42)$, we get
\begin{equation}
2a(s)\left \langle C,e_{1}\wedge e_{3}\right \rangle -b(s)\left \langle
C,e_{1}\wedge e_{2}\right \rangle =0
\end{equation}%
Combining $(38),(39)$ and $(42)$ we then have%
\begin{equation*}
f=4\left( a^{2}\left( s\right) -b^{2}\left( s\right) \right)
\end{equation*}%
that is, a smooth function $f$ depends only on $s.$ By differentiating $f$
with respect to $s$ and by using $(35)$ and $(27)$, we get%
\begin{equation}
f^{\prime }=2\varepsilon _{1}a(s)f
\end{equation}%
By differentiating $(38)$ with respect to $s$ and by using $%
(15),(27),(35),(36)$ and $(37)$ we have%
\begin{equation*}
a^{2}b=0
\end{equation*}%
or from $(35)$ we can write%
\begin{equation*}
\lambda a^{3}=0
\end{equation*}%
Since $a(s)\neq 0$, it follows that $\lambda =0.$ Then we obtain that $b=c=0.
$ Then the surface $M_{1}$ is a part of plane.

Thus we can give the following theorems.

\begin{theorem}
\label{teo1}Let $M_{1}$ be the flat rotation surface given by the
parametrization (7). Then $M_{1}$ has pointwise 1-type Gauss map if and only
if $M$ is either totally geodesic or parametrized by
\begin{equation}
\varphi \left( t,s\right) =\left(
\begin{array}{c}
\frac{1}{2}\varepsilon \left( \mu _{2}e^{\left( -b_{0}s+d\right) }-\mu
_{1}e^{-\left( -b_{0}s+d\right) }\right) \sinh t, \\
\frac{1}{2}\varepsilon \left( \mu _{2}e^{\left( -b_{0}s+d\right) }+\mu
_{1}e^{-(-b_{0}s+d)}\right) \cosh t, \\
\frac{1}{2}\varepsilon \left( \mu _{2}e^{\left( -b_{0}s+d\right) }+\mu
_{1}e^{-(-b_{0}s+d)}\right) \sinh t, \\
\frac{1}{2}\varepsilon \left( \mu _{2}e^{\left( -b_{0}s+d\right) }-\mu
_{1}e^{-\left( -b_{0}s+d\right) }\right) \cosh t%
\end{array}%
\right) \text{,\  \  \ }\mu _{1}\mu _{2}=-\frac{1}{b_{0}^{2}}
\end{equation}%
where $b_{0},$ $\mu _{1}$, $\mu _{2}$ and $d$ are real constants.
\end{theorem}

\begin{example}
\label{example1}Let $M_{1}$ be the flat rotation surface with pointwise
1-type Gauss map given by the parametrization (44). If we take as $b_{0}=-1,$
$\mu _{1}=-1,$ $\mu _{2}=1,$ $d=0$ and $\varepsilon =1$, then we obtain a
surface as follows:%
\begin{equation*}
\varphi \left( t,s\right) =\left( \cosh s\sinh t,\sinh s\cosh t,\sinh s\sinh
t,\cosh s\cosh t\right) .
\end{equation*}%
This surface is the product of two plane hyperbolas.
\end{example}

\begin{theorem}
\label{teo2}Let $M_{2}$ be the flat rotation surface given by the
parametrization (8). Then $M_{2}$ has pointwise 1-type Gauss map if and only
if $M_{2}$ is either totally geodesic or parametrized by
\begin{equation}
\text{ }\varphi \left( t,s\right) =\left(
\begin{array}{c}
\frac{1}{2}\varepsilon \left( \mu _{2}e^{\left( -b_{0}s+d\right) }+\mu
_{1}e^{-(-b_{0}s+d)}\right) \cos t, \\
\frac{1}{2}\varepsilon \left( \mu _{2}e^{\left( -b_{0}s+d\right) }+\mu
_{1}e^{-(-b_{0}s+d)}\right) \sin t, \\
\frac{1}{2}\varepsilon \left( \mu _{2}e^{\left( -b_{0}s+d\right) }-\mu
_{1}e^{-\left( -b_{0}s+d\right) }\right) \cos t, \\
\frac{1}{2}\varepsilon \left( \mu _{2}e^{\left( -b_{0}s+d\right) }-\mu
_{1}e^{-\left( -b_{0}s+d\right) }\right) \sin t%
\end{array}%
\right) \text{,\  \  \ }\mu _{1}\mu _{2}=-\frac{1}{b_{0}^{2}}
\end{equation}
\end{theorem}

\begin{example}
\label{example2}Let $M_{1}$ be the flat rotation surface with pointwise
1-type Gauss map given by the parametrization (44). If we take as $b_{0}=-1,$
$\mu _{1}=-1,$ $\mu _{2}=1,$ $d=0$ and $\varepsilon =1$, then we obtain a
surface as follows:%
\begin{equation*}
\varphi \left( t,s\right) =\varphi \left( t,s\right) =\left( \cosh s\cos
t,\cosh s\sin t,\cosh s\cos t,\cosh s\sin t\right) .
\end{equation*}%
This surface is the product of a plane circle and a plane hyperbola.
\end{example}

\begin{corollary}
\label{cor1}Let $M$\ be flat general rotation surface given by the
parametrization (7) or (8). If $M$ has pointwise 1-type Gauss map then the
Gauss map $G$ on $M$ is of 1-type.
\end{corollary}

\end{document}